\documentclass[11pt,leqno]{article}
\usepackage{amssymb}
\usepackage{amsmath}

\newtheorem{theorem}{Theorem}[section]

\numberwithin{equation}{section}

\newcommand{\A}{\mathcal{A}}
\newcommand{\f}{\mathcal{F}}
\newcommand{\g}{\mathcal{G}}

\begin{document}

\title{Decay in Time for a One Dimensional Two Component Plasma}
\author{Robert Glassey \\
Stephen Pankavich \\
Jack Schaeffer}
\date{\today}
\maketitle

\begin{center}
\end{center}

\begin{abstract}
The motion of a collisionless plasma is described by the
Vlasov-Poisson system, or in the presence of large velocities, the
relativistic Vlasov-Poisson system.  Both systems are considered in
one space and one momentum dimension, with two species of oppositely
charged particles. A new identity is derived for both systems and is
used to study the behavior of solutions for large times.
\end{abstract}

\section{Introduction}

Consider the Vlasov-Poisson system (which we shall abbreviate as
``VP''):

\begin{equation}\label{VP}
\left \{ \begin{gathered}
\partial_t f + v \ \partial_x f  + E(t,x) \ \partial_v f = 0,\\
\partial_t g + \frac{v}{m} \ \partial_x g  - E(t,x) \ \partial_v g = 0,\\
\rho(t,x) = \int \left ( f(t,x,v) - g(t,x,v) \right ) \ dv, \\
E(t,x) = \frac{1}{2} \left ( \int_{-\infty}^x \rho(t,y) \ dy -
\int_x^\infty \rho(t,y) \ dy \right ).
\end{gathered} \right.
\end{equation}
Here $t \geq 0$ is time, $x \in \mathbb{R}$ is position, $v \in
\mathbb{R}$ is momentum, $f$ is the number density in phase space of
particles with mass one and positive unit charge, and $g$ is the
number density of particles with mass $m > 0$ and negative unit
charge. The effect of collisions is neglected.  The initial
conditions $$ f(0,x,v) = f_0(x,v) \geq 0$$ and $$g(0,x,v) = g_0(x,v)
\geq 0$$ for $(x,v) \in \mathbb{R}^2$ are given where it is assumed
that $f_0,g_0 \in C^1(\mathbb{R}^2)$ are nonnegative, compactly
supported, and satisfy the neutrality condition
\begin{equation}
\label{neutrality} \iint f_0 \ dv dx = \iint g_0 \ dv\, dx.
\end{equation}
Using the notation
$$ \hat{v}_m = \frac{v}{\sqrt{m^2 + v^2}},$$ the relativistic
Vlasov-Poisson system (abbreviated ``RVP'') is
\begin{equation}
\label{RVP} \left \{ \begin{gathered}
\partial_t f + \hat{v}_1 \ \partial_x f  + E \ \partial_v f = 0,\\
\partial_t g + \hat{v}_m \ \partial_x g  - E \ \partial_v g = 0,\\
\rho(t,x) = \int \left ( f - g \right ) \ dv, \\
E(t,x) = \frac{1}{2} \left ( \int_{-\infty}^x \rho \ dy -
\int_x^\infty \rho \ dy \right ).
\end{gathered} \right.
\end{equation}

It is well known that solutions of (\ref{VP}) and (\ref{RVP}) remain
smooth for all $t \geq 0$ with $f(t,\cdot, \cdot)$ and
$g(t,\cdot,\cdot)$ compactly supported for all $t \geq 0$.  In fact,
this is known for the three-dimensional version of (\ref{VP})
(\cite{LP},\cite{Pfa}), but not for the three-dimensional version of
(\ref{RVP}).  The literature regarding large time behavior of
solutions is quite limited.  Some time decay is known for the
three-dimensional analogue of (\ref{VP})(\cite{GS}, \cite{IR},
\cite{Per}).  Also, there are time decay results for (\ref{VP}) (in
dimension one) when the plasma is monocharged (set $g \equiv 0$)
(\cite{BKR}, \cite{BFFM}, \cite{Sch}).  In the work that follows,
two species of particles with opposite charge are considered, thus
the methods used in \cite{BKR}, \cite{BFFM}, and \cite{Sch} do not
apply.  References \cite{DD}, \cite{Dol}, and \cite{DR} are also
mentioned since they deal with time-dependent rescalings and time
decay for other kinetic equations.\\

In the next section an identity is derived for (\ref{VP}) that shows
certain positive quantities are integrable in $t$ on the interval
$[0,\infty)$.  The identity is modified to address (\ref{RVP}) also,
but the results are weaker.  These identities seem to be linked to
the one-dimensional situation and do not readily generalize to
higher dimension.  Additionally, it is not clear if there is an
extension which allows for more than two species of particles.
 However, since this model allows for attractive forces between ions
of differing species, it is sensible to expect that additional
species of ions will only strengthen repulsive forces and cause
solutions to decay faster in time.  In Section $3$, the $L^4$
integrability of both the positive and negative charge is derived
and used to show time decay of the local charge. Finally, in Section
$4$, the main identity and $L^4$ integrability will be used to show
decay in time of the electric field for both (\ref{VP}) and
(\ref{RVP}).

\section{The Identity}

The basic identities for (\ref{VP}) and (\ref{RVP}) will be derived
in this section.  The following theorem lists their main
consequences :

\begin{theorem} \label{identity} Assume that $f_0$ and $g_0$ are nonnegative,
compactly supported, $C^1$, and satisfy (\ref{neutrality}).  Then,
for a solution of (\ref{VP}), there exists $C > 0$ depending only on
$f_0,g_0$, and $m$ such that $$ \int_0^\infty \iiint f(t,x,w)
f(t,x,v) (w-v)^2 \ dw\, dv\, dx\, dt \leq C,$$
$$ \int_0^\infty \iiint
g(t,x,w) g(t,x,v) (w-v)^2 \ dw\, dv\, dx\, dt \leq C,$$ and
$$\int_0^\infty \int E^2 \int (f + g) \ dv\, dx\,
dt \leq C.$$ For a solution of (\ref{RVP}) there is $C > 0$
depending only on $f_0,g_0$, and $m$ such that $$ \int_0^\infty
\iiint f(t,x,w) f(t,x,v) (w-v)(\hat{w}_1 - \hat{v}_1) \ dw\, dv\,
dx\, dt \leq C,$$
$$ \int_0^\infty \iiint g(t,x,w) g(t,x,v) (w-v)(\hat{w}_m - \hat{v}_m) \ dw\, dv\, dx\, dt \leq C,$$
and
$$\int_0^\infty \int E^2 \int (f + g) \ dv\, dx\,
dt \leq C.$$  Moreover, $(w-v)(\hat{w}_m - \hat{v}_m) \geq 0$ for
all $w,v \in \mathbb{R}$, $m > 0$.
\end{theorem}
{\bf Proof:} Suppose
\begin{equation}
\label{2.1}
\partial_t a + \omega(v) \partial_xa + B(t,x) \partial_v a = 0
\end{equation}
where $a(t,x,v), \omega(v), B(t,x)$ are $C^1$, and
$a(t,\cdot,\cdot)$ is compactly supported for each $t \geq 0$. Let
$$ A(t,x) = \int a \ dv$$ and $$\A(t,x) =
\int_{-\infty}^x A(t,y) \ dy.$$  Note that $\partial_x \A = A$ and
\begin{eqnarray*}
\partial_t \A & = &
-\int_{-\infty}^x \int \left ( \omega(v) \partial_y a(t,y,v) +
B(t,y) \partial_v a(t,y,v) \right ) \ dy\,dv \\
& = & - \int \omega(v) a(t,x,v) \ dv.
\end{eqnarray*}
By (\ref{2.1}) it follows that
\begin{equation}
\label{2.2}
\begin{split} 0 & = \A(t,x) \int v \left ( \partial_t a +
\omega(v) \partial_x a + B(t,x) \partial_v a \right ) \ dv \\
& =: I + II + III.
\end{split}
\end{equation}
Then,
\begin{equation}\label{2.3}
\begin{split}
I & = \partial_t\left (\A\int a v \ dv \right ) -
(\partial_t \A) \int a v \ dv \\
& = \partial_t \left (\A\int a v \ dv \right ) + \left (\int a v \
dv \right) \left (\int a \omega(v) \ dv \right ),
\end{split}
\end{equation}
\begin{equation}\label{2.4}
II = \partial_x \left ( \A \int a v \omega(v) \ dv \right ) - A
\int a v \omega(v) \ dv,
\end{equation}
and
\begin{equation}
\label{2.5}
\begin{split}
III &= -\A B \int a \ dv \\
& = -\A B A \\
& = -B \partial_x(\frac{1}{2} \A^2) \\
& = -\partial_x \left (\frac{1}{2} \A^2 B \right ) + \frac{1}{2}\A^2
\partial_x B.
\end{split}
\end{equation}
Using (\ref{2.3}), (\ref{2.4}), and (\ref{2.5}) in (\ref{2.2}) we
get

\begin{equation}
\label{2.6}
\begin{split}
0 &= \partial_t \left ( \A \int a v \ dv \right) +
\partial_x \left ( \A \int a v \omega(v) \ dv -
\frac{1}{2} \A^2 B \right ) \\
& \quad + \left ( \int a v \ dv \right ) \left (\int a \omega(v) \
dv \right ) - A \int a v \omega(v) \ dv + \frac{1}{2} \A^2
\partial_x B.
\end{split}
\end{equation}

Next consider (\ref{VP}) and let $$ F(t,x) := \int f \ dv, \quad
G(t,x) := \int g \ dv,$$ and $$ \f(t,x) := \int_{-\infty}^x F(t,y) \
dy, \quad \g := \int_{-\infty}^x G(t,y) \ dy.$$ Applying (\ref{2.6})
twice, once with $a = f$, $\omega(v) = v$ and $B = E$, and once with
$a = g$, $\omega(v) = \frac{v}{m}$, and $B = -E$, and adding the
results we find
\begin{equation}
\label{2.7}
\begin{split}
0 & = \partial_t \left ( \f \int f v \ dv + \g \int g v \ dv
\right ) + \partial_x \left ( \f \int f v^2 \ dv + m^{-1} \g \int
g v^2 \ dv \right ) \\
& \quad - \partial_x \left ( \frac{1}{2} \f^2 E - \frac{1}{2} \g^2
E \right ) + \left ( \int f v \ dv \right)^2 + m^{-1} \left (\int
g v \ dv \right )^2 \\
& \quad - F\int f v^2 \ dv - m^{-1}G \int g v^2 \ dv +
\frac{1}{2}\f^2\rho - \frac{1}{2}\g^2\rho.
\end{split}
\end{equation}
It follows directly from (\ref{VP}) and (\ref{neutrality}) that
$$\int \rho(t,x) \ dx = \int \rho(0,x) \ dx = 0$$ and hence that
$E \rightarrow 0$ as $\vert x \vert \rightarrow \infty$.  In
addition, $$E = \f - \g.$$  Hence
\begin{eqnarray*}
\int (\f^2 - \g^2) \rho \ dx & = & \int (\f + \g) E \ \partial_x E
\ dx \\
& = & -\frac{1}{2} \int \partial_x(\f + \g) E^2 \ dx \\
& = & -\frac{1}{2} \int(F + G) E^2 \ dx.
\end{eqnarray*}
Integration of (\ref{2.7}) in $x$ yields
\begin{equation}
\label{2.8}
\begin{split}
0 & = \frac{d}{dt} \left ( \int \f \int f v \ dv\,dx + \int \g \int
g v \ dv\,dx \right ) \\
& \quad + \int \left [ \left ( \int f v \ dv \right)^2  - F\int f
v^2 \ dv  + m^{-1} \left (\left (\int g v \ dv \right )^2 - G \int
g v^2 \ dv \right ) \right ] \ dx \\
& \quad - \frac{1}{4}\int (F+G) E^2 \ dx.
\end{split}
\end{equation}
Notice that exchanging $w$ and $v$ we can write
\begin{eqnarray*}
-\left (\int f v \ dv \right)^2 + F\int f v^2 \ dv & = & \left (
\int f(t,x,w) \ dw \right ) \left ( \int f(t,x,v) v^2 \ dv \right
) \\
& & \quad - \left ( \int f(t,x,w) w \ dw \right ) \left (\int
f(t,x,v) v \ dv \right ) \\
& = & \iint f(t,x,w) f(t,x,v) \left ( \frac{1}{2}w^2 +
\frac{1}{2}v^2 - wv \right ) \ dw\,dv \\
& = & \frac{1}{2} \iint f(t,x,w) f(t,x,v) (w - v)^2 \ dw\, dv
\end{eqnarray*}
and similarly for $g$.  Thus, (\ref{2.8}) yields
\begin{equation}
\label{2.9}
\begin{split}
\frac{d}{dt} \left ( \int \f \int f v \ dv\,dx + \int \g \int g v \
dv\,dx \right ) & = \frac{1}{2} \iiint f(t,x,w)
f(t,x,v) (w - v)^2 \ dw\, dv\, dx \\
& \quad + \frac{1}{2} m^{-1} \iiint g(t,x,w) g(t,x,v) (w - v)^2 \
dw\, dv\, dx \\
& \quad + \frac{1}{4} \int (F + G) E^2 \ dx \\
& \geq 0.
\end{split}
\end{equation}
Consider the energy $$ \iint (f + m^{-1}g) v^2 \ dv\,dx + \int E^2 \
dx.$$  Note that due to (\ref{neutrality}), $E(t,\cdot)$ is
compactly supported and $\int E^2 \ dx$ is finite (this would fail
without (\ref{neutrality})).  It is standard to show that the energy
is constant in $t$.  Similarly $\iint f \ dv\, dx = \iint g \ dv\,
dx$ is constant and $f, g \geq 0$.  Hence,
\begin{eqnarray*}
\left \vert \int \f \int f v \ dv\, dx \right \vert & \leq &
C\iint f \vert v \vert \ dv\, dx \\
& \leq & C \left ( \iint f \ dv\, dx \right )^\frac{1}{2} \ \left
( \iint f v^2 \ dv\, dx \right )^\frac{1}{2} \\
& \leq & C
\end{eqnarray*}
and similarly for $g$.  Now it follows from (\ref{2.9}) that
\begin{equation}
\label{2.10} \int_0^\infty \iiint \left ( f(t,x,w) f(t,x,v) +
g(t,x,w) g(t,x,v) \right ) (w - v)^2 \ dw\,dv\,dx\,dt \leq C
\end{equation}
and
\begin{equation}
\label{2.11} \int_0^\infty \int (F + G) E^2 \ dx \,dt \leq C.
\end{equation}

Next consider (\ref{RVP}).  Applying (\ref{2.6}) twice, once with $a
= f$, $\omega(v) = \hat{v}_1$, $B = E$ and once with $a = g$,
$\omega(v) = \hat{v}_m$, $B = -E$, and adding the results we get
\begin{eqnarray*}
0 & = & \partial_t \left ( \f \int f v \ dv + \g \int g v \ dv
\right ) + \partial_x \left ( \f \int f v \hat{v}_1 \ dv + \g \int g v \hat{v}_m \ dv \right ) \\
& & \quad - \partial_x \left ( \frac{1}{2} \f^2 E - \frac{1}{2} \g^2
E \right ) + \left ( \int f v \ dv \right) \left (\int f \hat{v}_1 \
dv \right ) + \left (\int g v \ dv \right ) \left ( \int g \hat{v}_m \ dv \right ) \\
& & \quad - F\int f v \hat{v}_1 \ dv - G \int g v \hat{v}_m \ dv +
\frac{1}{2}\f^2\rho - \frac{1}{2}\g^2\rho.
\end{eqnarray*}
Proceeding as before we obtain the result
\begin{equation}
\label{2.12}
\begin{split}
0 & = \frac{d}{dt} \left ( \int \f \int f v \
dv\,dx + \int \g \int g v \ dv\,dx \right ) \\
& \quad + \int \left [ \left (\int f v \ dv \right ) \left (\int f
\hat{v}_1 \ dv \right) - F\int f v \hat{v}_1 \ dv \right. \\
& \quad + \left. \left (\int g v \ dv \right ) \left (\int g
\hat{v}_m \ dv \right) - G \int g v
\hat{v}_m \ dv \right ] \ dx \\
& \quad - \frac{1}{4} \int (F + G) E^2 \ dx.
\end{split}
\end{equation}
Note that
\begin{equation*}
\begin{gathered}
-\left (\int f v \ dv \right ) \left (\int f \hat{v}_1 \ dv
\right) + F\int f v \hat{v}_1 \ dv \\
= \left (\int f(t,x,w) \ dw \right) \left (\int f(t,x,v) v
\hat{v}_1 \ dv \right ) - \left (\int f(t,x,w) w \ dw \right)
\left (\int f(t,x,v) \hat{v}_1 \ dv \right ) \\
= \frac{1}{2} \iint f(t,x,w) f(t,x,v) (v \hat{v}_1 + w
\hat{w}_1 - w \hat{v}_1 - v\hat{w}_1) \ dw\,dv \\
= \frac{1}{2} \iint f(t,x,w) f(t,x,v) (w - v) (\hat{w}_1 -
\hat{v}_1) \ dw\,dv.
\end{gathered}
\end{equation*}
By the mean value theorem for any $w$ and $v$, there is $\xi$
between them such that $$ \hat{w}_1 - \hat{v}_1 = (1 +
\xi^2)^{-\frac{3}{2}} (w - v)$$ and hence
\begin{equation}
\label{2.13} (w - v)(\hat{w}_1 - \hat{v}_1) = (1 +
\xi^2)^{-\frac{3}{2}} (w - v)^2 \geq 0.
\end{equation}
Similar results hold for $g$.  For solutions of (\ref{RVP}) $$ \iint
\left ( f \sqrt{1 + \vert v \vert^2} + g\sqrt{m^2 + \vert v \vert^2}
\right ) \ dv\, dx + \frac{1}{2} \int E^2 \ dx = const.$$ and mass
is conserved so
\begin{eqnarray*}
\left \vert \int \f \int f v \ dv\,dx + \int \g \int g v \ dv\,dx
\right \vert & \leq & C \iint f \vert v \vert \ dv\,dx + C \iint g
\vert v \vert \ dv\,dx \\
& \leq & C.
\end{eqnarray*}
Hence it follows by integrating (\ref{2.12}) in $t$ that
$$\int_0^\infty \iiint f(t,x,w) f(t,x,v) (w - v)(\hat{w}_1
- \hat{v}_1) \ dw\,dv\,dx\,dt \leq C,$$ $$ \int_0^\infty \iiint
g(t,x,w) g(t,x,v)(w - v)(\hat{w}_m - \hat{v}_m) \ dw\,dv\,dx\,dt
\leq C,$$ and $$\int_0^\infty \int (F + G)E^2 \ dx\,dt \leq C.$$
Theorem \ref{identity} now follows.
\begin{flushright} $\square$
\end{flushright}

\section{Decay Estimates}

In this section we will derive some consequences of the identity
from the previous section.  We begin by taking $m=1$ and defining
$\hat{v} := \hat{v}_m = \hat{v}_1$.  Consider solutions to either  the
system (\ref{VP}) or the system (\ref{RVP}), and define as above
$$F(t,x)=\int f(t,x,v)\,dv, \quad  G(t,x)=\int g(t,x,v)\,dv.$$

\begin{theorem}\label{L4bound}
Let $f,\,g$ satisfy the VP system (\ref{VP}).  Assume that the
data functions $f_0,\,g_0$ satisfy the hypotheses of Theorem
\ref{identity}. Then
$$\int_0^{\infty}\int F^4(t,x)\,dx\,dt<\infty$$
and $$\int_0^{\infty}\int G^4(t,x)\,dx\,dt<\infty.$$ When
$f,\,g$ satisfy the RVP system (\ref{RVP}) and the data
functions $f_0,\,g_0$ satisfy the hypotheses of Theorem
\ref{identity} we have
$$\int_0^\infty \left (\int F(t,x)^\frac{7}{4} \ dx
\right )^4 \ dt < \infty$$
with the same result valid for $G$.
\end{theorem}
{\bf Proof:} Consider the classical case (\ref{VP}).  By Theorem \ref{identity} we know that
$$k(t,x) := \iint  (w-v)^2f(t,x,v) f(t,x,w)\,dv\,dw$$
is integrable over all $x,\,t$.  Next we partition the set of integration:
$$F(t,x)^2 = \iint f(t,x,v) f(t,x,w)\,dv\,dw=\int_{|v-w|<R}+\int_{|v-w|>R} =: I_1+I_2.$$
Clearly we have $I_2\le R^{-2}k(t,x)$.  In the integral for $I_1$
we write
$$\int_{|v-w|<R}f(t,x,w)\,dw=\int_{v-R}^{v+R}f(t,x,w)\,dw\le 2\|f_0\|_{\infty}R.$$
Thus
$$I_1\le c\cdot R \cdot F.$$
Set $RF=R^{-2}k$ or $R=k^{1\over 3}F^{-{1\over 3}}$.  Then
$F^4(t,x)\le ck(t,x)$ so $F^4(t,x)$ is integrable over all
$x,\,t$. The result for $G$ is exactly the same.

Now we will find by a similar process the corresponding estimate for
solutions to the relativistic version (\ref{RVP}). To derive it we
will use the estimate from (\ref{2.13}), which implies that for
$1+|v|+|w|\le S$, there is a constant $c>0$ such that
 $$(v-w)(\hat v - \hat w)\ge cS^{-3}|v-w|^2.$$
From Theorem \ref{identity} with $m=1$ we know that
$$k_r(t,x) := \iint  (v-w)(\hat v - \hat w)f(t,x,v) f(t,x,w)\,dv\,dw$$
is integrable over all $x,\,t$.  Now write
$$F(t,x)^2=\iint f(t,x,v) f(t,x,w)\,dv\,dw=\int_{(v-w)(\hat v - \hat w)<R}
+\int_{(v-w)(\hat v - \hat w)>R} =: I_1+I_2.$$ Clearly  $I_2\le
R^{-1}k_r(t,x)$.  To estimate $I_1$ we partition it as
$$I_1=\iint_{(v-w)(\hat v - \hat w)<R \atop
1+|v|+|w|<S} f(t,x,v) f(t,x, w)\,dv\,dw+ \iint_{(v-w)(\hat v -
\hat w)<R\atop 1+|v|+|w|>S} f(t,x,v) f(t,x,w)\,dv\,dw =: I_1' +
I_1''.$$ On $I_1'$ we have by the above estimate
$$R\ge (v-w)(\hat v - \hat w)\ge c|v-w|^2S^{-3}.$$
Therefore on $I_1'$ we have $|v-w|<cR^{1\over 2}S^{3\over 2}$
so that
$$I_1'\le c\int f(t,x,v)  \int_{v-cR^{1\over 2}S^{3\over 2}}
^{v+cR^{1\over 2}S^{3\over 2}}f(t,x,w)\,dw\,dv\le c\cdot
F(t,x)\cdot R^{1\over 2}S^{3\over 2}.$$ $I_1''$ is more
troublesome.  By the energy and mass bounds,
$$I_1''\le S^{-1}\iint (1+|v|+|w|)f(t,x,v) f(t,x,w)\,dv\,dw\le cS^{-1}e(t,x)
F(t,x)$$
where $e(t,x)=\int \sqrt{1+v^2}f(t,x,v)\,dv$.
Find $S$ first by setting
$$F(t,x)\cdot R^{1\over 2}S^{3\over 2}=S^{-1}e(t,x)F(t,x),$$
that is,
$$S=e(t,x)^{2\over 5}R^{-{1\over 5}}.$$
Thus we get for $I_1$ the bound
$$I_1\le cS^{-1}e(t,x)F(t,x)= cF(t,x)R^{1\over 5}e(t,x)^{3\over 5}.$$
Above we had $I_2\le R^{-1}k_r(t,x)$.  So now set
$$F(t,x)R^{1\over 5}e(t,x)^{3\over 5}=R^{-1}k_r(t,x)$$
to find $R$.  The result is
$$R=k_r(t,x)^{5\over 6}F(t,x)^{-{5\over 6}}e(t,x)^{-{1\over 2}}.$$
Finally then
$$F(t,x)^2\le cR^{-1}k_r(t,x)=ck_r(t,x)^{1\over 6}F(t,x)^{5\over 6}e(t,x)^{1\over 2}$$
which is the same as ${F(t,x)^{7}\over e(t,x)^3}\le ck_r(t,x)$.  At
this point we may integrate in time to produce the result
\begin{equation}
\label{L4RVP1} \int_0^{\infty}\int{(\int f(t,x,v)\,dv)^{7}\over
(\int \sqrt{1+v^2}f(t,x,v)\,dv)^3}\,dx\,dt<\infty.
\end{equation}
Alternatively we can isolate $F(t,x)^7$ on the left side to find
$F(t,x)^7 \leq c k_r(t,x) e(t,x)^3$.  Then raise both sides to the
$\frac{1}{4}$th power, integrate in $x$ and use H\"older's inequality.
Hence we get the bound
$$ \int F(t,x)^\frac{7}{4} \ dx \leq \left (\int k_r(t,x) \
dx \right )^\frac{1}{4} \left (\int e(t,x) \ dx \right
)^\frac{3}{4}.$$  We use the time--independent bound on $\int e(t,x)
\ dx$ from conservation of energy to get the estimate
$$ \int F(t,x)^\frac{7}{4} \ dx \leq C \left ( \int k_r(t,x) \ dx
\right )^\frac{1}{4}.$$  Finally, we raise both sides to the $4$th
power and integrate in time to produce the result
\begin{equation}
\label{L4RVP2} \int_0^\infty \left (\int F(t,x)^\frac{7}{4} \ dx
\right )^4 \ dt < \infty.
\end{equation}
This is corresponding estimate for solutions to (\ref{RVP}).
\begin{flushright} $\square$
\end{flushright}
Now we will  use these estimates
to show that the local charges tend to $0$ as
$t\to \infty$ for solutions to both sets of equations.
\bigskip\noindent
\begin{theorem} \label{Localcharge}
Let $f,\,g$ be solutions to either the classical VP system (\ref{VP})
or to the relativistic RVP system (\ref{RVP}) for
which the assumptions of Theorem 3.1 hold.  Then for any fixed $R
> 0$ the local charges satisfy
$$\lim_{t\to \infty}\int_{|x|<R}F(t,x)\,dx=
\lim_{t\to \infty}\int_{|x|<R}G(t,x)\,dx=0.$$
\end{theorem}
{\bf Proof:} We begin with solutions to the classical equation (\ref{VP}).
From above we know that $$\int_0^{\infty}\int
F^4(t,x)\,dx\,dt<\infty.$$ By the H\"older inequality
$$\int_{|x|<R}F(t,x)\,dx\le \left(\int
F^4(t,x)\,dx\right)^{1/4}(2R)^{3/4}$$ and therefore
\begin{equation}\label{holder4}
\int_0^{\infty}\left[\int_{|x|<R}F(t,x)\,dx\right]^4\,dt<\infty.
\end{equation}
Now by the Vlasov equation for $f$
$$F_t=-\int (vf_x+Ef_v)\,dv = -\partial_x \int vf\,dv.$$
Integrate this  formula in $x$ over $|x|<R$:
$$\partial_t \int_{|x|<R}F(t,x)\,dx=- \int_{|x|<R}\partial_x \int vf\,dv\,dx=
-\int vf(t,R,v)\,dv+\int vf(t,-R,v)\,dv.$$
Call
$$j_f(t,x)=\int vf(t,x,v)\,dv.$$
Then $j_f(t,x)$ is boundedly integrable over all $x$ by the
energy. Next we compute
\begin{eqnarray*}
\partial_t \left[\int_{|x|<R}F(t,x)\,dx\right]^4 & = &
4\left[\int_{|x|<R}F(t,x)\,dx\right]^3\int_{|x|<R}F_t(t,x)\,dx\\
& = &
4\left[\int_{|x|<R}F(t,x)\,dx\right]^3\left[-j_f(t,R)+j_f(t,-R)\right].
\end{eqnarray*}
For $0<R_1<R_2$ integrate this in $R$ over $[R_1,R_2]$:
\begin{eqnarray*}
\frac{d}{dt} \int_{R_1}^{R_2}\left[\int_{|x|<R}F(t,x)\,dx\right]^4\,
dR & = & 4\int_{R_1}^{R_2}\left[\int_{|x|<R}F(t,x)\,dx\right]^3
\left[-j_f(t,R)+j_f(t,-R)\right]\,dR \\
&\leq & 4\left[\int_{|x|<R_2}F(t,x)\,dx\right]^3\int_{R_1}^{R_2}
\left|-j_f(t,R)+j_f(t,-R)\right|\,dR \\
&\leq & c\left[\int_{|x|<R_2}F(t,x)\,dx\right]^3
\end{eqnarray*}
for some constant $c$ depending only on the data.  For
$t_2>t_1>1$ multiply this by $t-t_1$ and integrate in $t$ over
$[t_1,t_2]$:
$$\int_{t_1}^{t_2}(t-t_1)\partial_t \int_{R_1}^{R_2}\left[\int_{|x|<R}F(t,x)\,dx\right]^4\,dR\,dt
\le
c\int_{t_1}^{t_2}(t-t_1)\left[\int_{|x|<R_2}F(t,x)\,dx\right]^3\,dt.
$$
Integrating the left side by parts we get
$$(t_2-t_1) \int_{R_1}^{R_2}\left[\int_{|x|<R}F(t_2,x)\,dx\right]^4\,dR
-\int_{t_1}^{t_2}\int_{R_1}^{R_2}\left[\int_{|x|<R}F(t,x)\,dx\right]^4\,dR\,dt.$$
Now take $t_2=t,\,t_1=t-1$.  Then we have
\begin{equation}\label{L4bound2} \begin{array}{rcl}
\displaystyle
\int_{R_1}^{R_2}\left[\int_{|x|<R}F(t,x)\,dx\right]^4\,dR & \leq & \displaystyle \int_{t-1}^{t}\!\int_{R_1}^{R_2}\left[\int_{|x|<R}F(t,x)\,dx\right]^4\,dR\,dt\\
& \ & \displaystyle \quad + c\int_{t-1}^{t}(t-t_1)\left[\!\int_{|x|<R_2}F(t,x)\,dx\right]^3\,dt \\
&\leq & \displaystyle \int_{t-1}^{t}(R_2-R_1)\left[\int_{|x|<R_2}F(t,x)\,dx\right]^4\,dt \\
& \ & \displaystyle \quad +
c\int_{t-1}^{t}\left[\int_{|x|<R_2}F(t,x)\,dx\right]^3\,dt.
\end{array} \end{equation}
Now take $R_2=2R_1 = 2R$ say.  Then the left side of
(\ref{L4bound2}) is bounded below by
$$(R_2-R_1) \left[\int_{|x|<R_1}F(t,x)\,dx\right]^4=R
\left[\int_{|x|<R}F(t,x)\,dx\right]^4$$ and we claim that the right
side tends to 0 as $t\to \infty$.  This is clear for the first term
on the right of (\ref{L4bound2}) by use of (\ref{holder4}). The
second term goes to 0, as well, since
$$\int_{t-1}^{t}\left[\int_{|x|<R_2}F(t,x)\,dx\right]^3\,dt \leq
\left(\int_{t-1}^{t}\left[\int_{|x|<R_2}F(t,x)\,dx\right]^4\,dt\right)^{3/4}
\cdot \left(\int_{t-1}^t\,dt\right)^{1/4}.$$ The same computation
establishes the estimate for $G$, and the result now follows in the classical
case.

The proof for the relativistic case is similar.
From above we know that
$$\label{L4RVP2} \int_0^\infty \left (\int F(t,x)^\frac{7}{4} \ dx
\right )^4 \ dt < \infty.$$
By the H\"older inequality
$$\int_{|x|<R}F(t,x)\,dx\le c_R\left(\int
F^\frac{7}{4}(t,x)\,dx\right)^\frac{4}{7}$$
and therefore
$$\int_0^{\infty}\left[\int_{|x|<R}F(t,x)\,dx\right]^7\,dt<\infty.$$
Using the Vlasov equation (\ref{RVP}) for $f$ we have
$$F_t=-\int (\hat vf_x+Ef_v)\,dv = -\partial_x \int \hat vf\,dv.$$
Integrate this  formula in $x$ over $|x|<R$:
$$\partial_t \int_{|x|<R}F(t,x)\,dx=- \int_{|x|<R}\partial_x \int \hat vf\,dv\,dx=
-\int \hat vf(t,R,v)\,dv+\int \hat vf(t,-R,v)\,dv.$$
Call
$$j_f^r(t,x)=\int \hat vf(t,x,v)\,dv.$$
Then $j_f^r(t,x)$ is boundedly integrable over all $x$ by the
mass bound. Next we compute
\begin{eqnarray*}
\partial_t \left[\int_{|x|<R}F(t,x)\,dx\right]^7 & = &
7\left[\int_{|x|<R}F(t,x)\,dx\right]^6\int_{|x|<R}F_t(t,x)\,dx\\
& = &
7\left[\int_{|x|<R}F(t,x)\,dx\right]^6\left[-j_f^r(t,R)+j_f^r(t,-R)\right].
\end{eqnarray*}
The proof now concludes  exactly as in the classical case.
\begin{flushright} $\square$
\end{flushright}

\section{Time Decay of Electric Field}
We conclude the paper with results concerning the time
integrability and decay of the electric field for both the
classical and relativistic systems, (\ref{VP}) and (\ref{RVP}).

\begin{theorem} \label{E3int}
Let the assumptions of Theorem \ref{L4bound} hold and consider
solutions $f,\,g$ to either (\ref{VP}) or (\ref{RVP}). Then
$$\int_0^{\infty}\Vert E(t) \Vert_{\infty}^3\,dt<\infty.$$
\end{theorem}
{\bf Proof:}  This will follow immediately from the result in
Theorem \ref{identity} that
$$Q(t) := \int_{-\infty}^{\infty} E^2(t,x)\left[F(t,x)+G(t,x)\right]\,dx$$
is integrable in time.  Indeed by the equation $E_x=\rho =
\int(f-g)\,dv=F-G$, we have
$${\partial\over \partial x}E^3=3E^2\rho=3E^2(F-G).$$
Integrate in $x$ to get
$$E^3(t,x)=\int_{-\infty}^x 3E^2(F-G)\,dx$$
so that
\begin{equation}
\label{E3int} |E(t,x)|^3\le \int_{-\infty}^{\infty} 3E^2(F+G)\,dx
= 3Q(t)
\end{equation}
and the result follows as claimed.
\begin{flushright}
$ \square $
\end{flushright}

Our final results will show that for solutions to the classical VP
system (\ref{VP}) and RVP system (\ref{RVP}), the electric field $E$
tends to 0 in the maximum norm.
\begin{theorem} \label{EdecayVP}
Let the assumptions of Theorem \ref{L4bound} hold and consider
solutions $f,\,g$ to the classical VP system (\ref{VP}).   Then
$$\lim_{t\to \infty}\|E(t)\|_{\infty}=0.$$
\end{theorem}
{\bf Proof:}  We will show that
$$\lim_{t\to \infty}Q(t)=0.$$
The conclusion will then follow from (\ref{E3int}). Since $Q(t)$
is integrable over $[0,\infty)$, $\liminf Q(t)=0.$ Therefore,
there is a sequence $t_n$ tending to infinity such that  $Q(t_n)
\to 0$ as $n \to \infty$.  As above, we denote
$$F(t,x)=\int f(t,x,v)\,dv, \quad  G(t,x)=\int g(t,x,v)\,dv.$$
Using $E_x=\rho=F-G$ and $E_t=-j=-\int v(f-g)\,dv$ we first
compute
\begin{eqnarray*}
\frac{dQ}{dt}& = & -2\int jE(F+G)\,dx + \int
E^2\partial_t(F+G)\,dx\\
& = & -2\int jE(F+G)\,dx - \int E^2\partial_x\int v(f+g)\,dv
\,dx\\
& = & -2\int jE(F+G)\,dx + 2\int \rho E\int v(f+g)\,dv\,dx.
\end{eqnarray*}
Now, $E$ is uniformly bounded because by definition in (\ref{VP}),
$$|E(t,x)|\le \int_{-\infty}^x (F+G)(t,x)\,dx\le
\int_{-\infty}^{\infty} (F+G)(t,x)\,dx\le \hbox{const.}$$ where the
last inequality follows by conservation of mass.   Therefore
$$\left|{dQ\over dt}\right| \le  c \int (F+G)\int |v|(f+g)\,dv \,dx.$$
Define $e$ to be the kinetic energy density, $$e(t,x) := \int v^2
(f+g)\,dv.$$ Then in the usual manner we get
\begin{eqnarray*}
\int |v|(f+g)dv & = & \int_{|v|<R}|v|(f+g)\,dv +
\int_{|v|>R}|v|(f+g)\,dv\\
&\leq & \Vert f+g \Vert_{\infty}\cdot R^2+R^{-1}e\\
&\leq & c(R^2+R^{-1}e).
\end{eqnarray*}
Choosing $R^3=e$ we find that
$$\int |v|(f+g)dv \le  c e^{2\over 3}(t,x)$$ and therefore
\begin{equation} \label{dQdt}
\left|{dQ\over dt}\right|  \le  c \int (F+G)e^\frac{2}{3}\,dx \le
c \Big(\int (F+G)^3 \,dx\Big)^{1\over 3}
\end{equation} by the H\"older inequality and the bound on
kinetic energy from Section $2$.  By interpolation, for suitable
functions $w$,
$$ \|w\|_3 \le   \|w\|_1^{\theta} \cdot \|w\|_4^{1-\theta}$$
where
$${1\over 3}={\theta \over 1} + {1-\theta \over 4}.$$
Therefore $\theta = {1\over 9}$.  Apply this to $w=F+G$ and use the boundedness of
$F+G$ in $L^1$ to get
$$ \|F+G\|_3 \le c \|F+G\|_4^{8\over 9}.  $$
Using this above we conclude that
$$\left|{dQ\over dt}\right| \le c\|F+G\|_4^{8\over 9}.$$
From Theorem \ref{L4bound} we know that $\int (F^4 + G^4)\,dx$ is
integrable in time. Thus $\left|{dQ\over dt}\right|^{9\over 2}$ is
integrable in time.  Now, for any $0  < R_1 < R_2$ write
$$Q(R_2)^{16\over 9} - Q(R_1)^{16\over 9}={16\over 9}\int_{R_1}^{R_2}Q(t)^{7\over 9}
\dot Q(t)\,dt.$$ By the H\"older inequality again, with $p={9\over
7}$ and $q={9\over 2}$,
$$\Big|Q(R_2)^{16\over 9} - Q(R_1)^{16\over 9}\Big| \le c\Big(\int_{R_1}^{R_2}Q(t)\,dt\Big)^{7/9} \cdot
\Big(\int_{R_1}^{R_2}|\dot Q(t)|^{9\over 2}\,dt\Big)^{2\over 9} \to 0
$$
as $R_1, \, R_2 \to \infty.$ Therefore the limit
$$\lim_{R \to \infty} Q(R)^{16\over 9}$$ exists and equals $\omega$, say.
By taking $R=t_n$  and letting $n \to \infty$ we get $\omega =0$.  This concludes
the proof.
\begin{flushright}
$\square$
\end{flushright}

\begin{theorem} \label{EdecayRVP}  Let the assumptions of Theorem 2.1 hold and
consider solutions $f,\,g$ to the relativistic VP system
(\ref{RVP}). Then, also in this case
$$\lim_{t\to \infty}\|E(t)\|_{\infty}=0.$$
\end{theorem}
{\bf Proof:} As is to be expected, the proof is similar to that of
Theorem \ref{EdecayVP}. From Theorem 2.1 we have again that $Q(t)$
is integrable in time, where exactly as in the non-relativistic case
$$Q(t)= \int_{-\infty}^{\infty} E^2(t,x)\left[F(t,x)+G(t,x)\right]\,dx.$$
In this situation we have $\rho=\int(f-g)\,dv$ and (with $m=1$)
$j=\int \hat v(f-g)\,dv$ where $\hat v={v\over \sqrt{1+v^2}}$ so
that $|\hat v|<1$.  The computation of the derivative in time is
now
\begin{eqnarray*}
\frac{dQ}{dt} &=& -2\int jE(F+G)\,dx + \int
E^2\partial_t(F+G)\,dx\\
& = & -2\int jE(F+G)\,dx - \int E^2\partial_x\int \hat v(f+g)\,dv
\,dx\\
& = & -2\int jE(F+G)\,dx + 2\int \rho E\int \hat v(f+g)\,dv\,dx.
\end{eqnarray*}
It follows that
$$\left|{dQ\over dt}\right| \le c\int |E|(F+G)^2\,dx\le c\int (F+G)^2\,dx$$
because $E$ is uniformly bounded. Call $e$ the relativistic kinetic
energy density,  $$e(t,x) = \int \sqrt{1+v^2} (f+g)\,dv.$$ Then as
above
\begin{eqnarray*}
F+G & = & \int (f+g)\,dv\\
& = & \int_{|v|<R}(f+g)\,dv + \int_{|v|>R}(f+g)\,dv\\
& \leq &  \Vert f+g \Vert_{\infty}\cdot 2R+R^{-1}e\\
&\leq & c(R+R^{-1}e).
\end{eqnarray*}
Hence with $R^2=e$ we find
that $F+G\le ce^{1\over 2}$.  Thus we see that
$$\left|{dQ\over dt}\right|  \le  c \int e\,dx \le c.$$
In view of Remark 1 below then, $Q(t)\to 0$ as $t\to \infty$ which
implies the result for $E$ as in the classical case.
\begin{flushright}
$ \square $
\end{flushright}

\noindent \underline{\textbf{Remarks:}}

\noindent \textbf{1.} Once $Q(t)$ is integrable in time, the
uniform boundedness of $\left|{dQ\over dt}\right|$ also implies
that $Q(t) \to  0$ as $t \to \infty$.  The estimate (\ref{dQdt})
provides the desired bound in the classical case because $(F+G)^3$
is dominated by the energy integral in this situation.\\

\noindent \textbf{2.} For solutions to (\ref{VP}) or (\ref{RVP}),
using interpolation with Theorem \ref{EdecayVP} or Theorem
\ref{EdecayRVP}, and the bound on $\Vert E(t) \Vert_2$ from energy
conservation we find $$ \lim_{t \rightarrow \infty} \Vert E(t)
\Vert_p = 0$$ for any $p > 2$.\\

\noindent \textbf{3.} We have been unable to find a rate of decay
for the maximum norm of $E$. For solutions to the classical
Vlasov--Poisson system in three space dimensions such a rate follows
from differentiating in time an expression essentially of the form
$$\iint |x-tv|^2(f+g)\,dv\,dx$$
(cf. \cite{IR}, \cite{Per}).  This estimate fails to imply time
decay in the current one--dimensional case.\\

\noindent \textbf{4.} An identity similar to that in the proof of
Theorem \ref{identity} holds for solutions to the ``one and
one--half--dimensional'' Vlasov--Maxwell system.  However we have
been unable to show that certain terms arising from the linear parts
of the differential operators have the proper
sign.\\

\noindent \textbf{5.} As stated in the introduction, such decay
theorems should be true for several species under the hypothesis of
neutrality. However we have been unable to achieve this
generalization for more than two species.\\

\noindent \textbf{6.} After suitable approximation, these results can be
seen to be valid for weak solutions as well.

\end{document}